\documentclass{article}
\usepackage{graphicx}
\include{graphics}
\usepackage{amscd,amsfonts,euscript,amsmath,amsxtra,amssymb,amsthm}
\usepackage[russian, english]{babel}
\newtheorem{theorem}{Theorem}[section]

\newtheorem{proposition}[theorem]{Proposition}

\theoremstyle{definition}
\newtheorem{definition}[theorem]{Definition}

\theoremstyle{remark}
\newtheorem{remark}[theorem]{Remark}

\newtheorem{convention}[theorem]{Convention}

\numberwithin{equation}{section}

\usepackage[all]{xy}

\def\E{{{\mathbb E }}}

\def\L{{{\mathbb L }}}

\def\II{{{\mathcal I}}}

\def\LL{{{\mathcal L}}}

\def\OO{{{\mathcal O}}}

\def\SS{{{\mathcal S}}}

\def\gl.dim{{{\rm gl.dim \,}}}

\def\ker{{{\rm ker \,}}}

\def\Pic{{{\rm Pic\,}}}
\begin{document}{\sl MSC~14J60, 14D20, 14M27}

{\sl UDC~512.722+512.723}
\medskip
\begin{center}
{\Large\sc On some isomorphism of compactifications of moduli scheme
of vector bundles} 
\end{center}
\medskip
\begin{center}
N.V.~Timofeeva

\smallskip

Yaroslavl' State University

Sovetskaya str. 14, 150000 Yaroslavl', Russia

e-mail: {\it ntimofeeva@list.ru}
\end{center}
\bigskip

\begin{quote}
{\it\qquad A morphism of the reduced Gieseker -- Maruyama moduli
functor
 (of semistable coherent torsion-free sheaves) to the reduced
moduli functor   of admissible semistable pairs with the same
Hilbert polynomial, is constructed. It is shown that main components
of reduced moduli scheme for semistable admissible pairs
$((\widetilde S, \widetilde L), \widetilde E)$ are isomorphic to
main components of reduced Gieseker -- Maruyama moduli scheme.


Bibliography: 9 items.}
\end{quote}

\section*{Introduction}
In the present article $S$ is a smooth irreducible projective
algebraic surface over an algebraically closed field $k$ of
characteristic  0, ${\mathcal O}_S$ is its structure sheaf, $E$
coherent torsion-free ${\mathcal O}_S$-module,
$E^{\vee}:={{\mathcal H}om}_{{\mathcal O}_S}(E, {\mathcal O}_S)$
dual ${\mathcal O}_S$-module. In this case  $E^{\vee}$ is
reflexive and, consequently, locally free. In the sequel we make
no difference between locally free sheaf and the corresponding
vector bundle, and both terms are used as synonyms. Let  $L$ be
very ample invertible sheaf on  $S$; it is fixed and will be
referred to as the polarization. The symbol  $\chi(\cdot)$ denotes
the Euler characteristic, $c_i(\cdot)$  $i$-th Chern class. Also
if $Y\subset X$ be the locally closed subscheme of the scheme $X$,
then $\overline Y$ be its scheme-theoretic closure in $X$.

\begin{definition} \cite{Tim3, Tim4} A polarized algebraic
scheme  $(\widetilde S, \widetilde L)$ is called {\it admissible} if
the scheme  $(\widetilde S,\widetilde L)$ satisfies one of following
two conditions

i) $(\widetilde S, \widetilde L) \cong (S,L)$,

ii) $\widetilde S \cong {\rm Proj \,} \bigoplus_{s\ge
0}(I[t]+(t))^s/(t^{s+1})$, where $I={{\mathcal F}itt}^0 {{\mathcal
E}xt}^2(\varkappa, {\mathcal O}_S)$ for Artinian quotient sheaf
$q: \bigoplus^r {\mathcal O}_S\twoheadrightarrow \varkappa$ of
length $l(\varkappa)\le c_2$, and $\widetilde L = L \otimes
(\sigma ^{-1} I \cdot {\mathcal O}_{\widetilde S})$ is very ample
invertible sheaf on the scheme $\widetilde S$. This polarization
$\widetilde L$ will be referred to as the {\it distinguished
polarization}.
\end{definition}
\begin{remark} In the further considerations, if necessary, we
replace  $L$ by its big enough tensor power. As shown in
\cite{Tim4}, this power can be chosen constant and fixed. All
Hilbert polynomials are compute with respect to these new  $L$ and
$\widetilde L$ correspondingly.
\end{remark}

\begin{definition} \cite{Tim4,Tim4E} $S$-{\it (semi)stable pair} $((\widetilde
S,\widetilde L), \widetilde E)$ is the following data:
\begin{itemize}
\item{$\widetilde S=\bigcup_{i\ge 0} \widetilde S_i$ admissible scheme, $\sigma: \widetilde S \to S$ morphism
which is called {\it canonical}, $\sigma_i: \widetilde S_i \to S$
its restrictions onto components $\widetilde S_i$, $i\ge 0;$}
\item{$\widetilde E$ vector bundle on the scheme
$\widetilde S$;}
\item{$\widetilde L \in Pic\, \widetilde S$ distinguished polarization;}
\end{itemize}
such that
\begin{itemize}
\item{$\chi (\widetilde E \otimes \widetilde
L^m)=rp(m),$ the polynomial  $p(m)$ and the rank  $r$ of the sheaf
 $\widetilde E$ are fixed;}
\item{ the sheaf
$\widetilde E$ on the scheme $\widetilde S$ is {\it stable
}(respectively, {\it semistable}) {\it in the sense of Gieseker,} то
i.e. for any proper subsheaf $\widetilde F \subset \widetilde E$ for
$m\gg 0$
\begin{eqnarray*}
\frac{h^0(\widetilde F\otimes \widetilde L^m)}{{\rm rank\,} F}&<&
\frac{h^0(\widetilde E\otimes \widetilde L^m)}{{\rm rank\,} E},
\\ (\mbox{\rm respectively,} \;\;
\frac{h^0(\widetilde F\otimes \widetilde L^m)}{{\rm rank\,} F}&\leq&
\frac{h^0(\widetilde E\otimes \widetilde L^m)}{{\rm rank\,} E}\;);
\end{eqnarray*}}
\item{for each additional component $\widetilde S_i, i>0,$
the sheaf $\widetilde E_i:=\widetilde E|_{\widetilde S_i}$ is {\it
quasi-ideal,} namely it has a description вида
\begin{equation}\label{quasiideal}\widetilde E_i=\sigma_i^{\ast}\ker q_0/tor\!s.\end{equation} для некоторого $q_0\in
\bigsqcup_{l\le c_2} {\rm Quot\,}^l \bigoplus^r {\mathcal O}_S$.
}\end{itemize}
\end{definition}

In the series of papers of the author \cite{Tim0}
--- \cite{Tim4E} the projective algebraic scheme $\widetilde M$ is
constructed as reduced moduli scheme for $S$-semistable pairs.

The scheme $\widetilde M$ contains an open subscheme $\widetilde
M_0$ which is isomorphic to the subscheme $M_0$ of
Gieseker-semistable vector bundles, in Gieseker -- Maruyama moduli
scheme $\overline M$ for semistable torsion-free sheaves with
Hilbert polynomial equal to $\chi(E \otimes L^m)=rp(m)$. We make
use of the Gieseker's definition of semistability.

\begin{definition}\cite{Gies} Coherent  ${\mathcal O}_S$-sheaf $E$ is {\it
stable} (respectively, {\it semistable})  if for any proper subsheaf
 $F\subset E$ of rank  $r'={\rm
rank\,} F$ for  $m\gg 0$ the following holds:
$$
\frac{\chi(E \otimes L^{m})}{r}>\frac{\chi(F\otimes L^m)}{r'},\;
{\mbox{\LARGE (}}{\mbox{\rm respectively,}} \; \frac{\chi(E \otimes
L^{m})}{r}\ge \frac{\chi(F\otimes L^m)}{r'}{\mbox{\LARGE )}}.
$$
\end{definition}

Let $E$ be a semistable locally free coherent sheaf. Then,
obviously, the sheaf  $I={{\mathcal F}itt}^0 {{\mathcal E}xt}^1(E,
{\mathcal O}_S)$ is trivial and  $\widetilde S \cong S$.
Consequently, $((\widetilde S, \widetilde L), \widetilde E) \cong
((S, L),E)$ and there is a bijection $\widetilde M_0 \cong M_0$.

Let  $E$ be a semistable nonlocally free sheaf, then the scheme
$\widetilde S$ contains reduced irreducible component $\widetilde
S_0$ such that  $\sigma_0:=\sigma|_{\widetilde S_0}: \widetilde S_0
\to S$ is a morphism of blowing up of the scheme  $S$ in the sheaf
of ideals  $I= {{\mathcal F}itt}^0 {{\mathcal E}xt}^1(E, {\mathcal
O}_S).$ Formation of the sheaf  $I$ is a way to characterize
singularities of the sheaf $E$, i.e. its difference from local
freeness. Indeed, the quotient sheaf $\varkappa:= E^{\vee \vee}/ E$
is an Artinian sheaf and its length is not greater then $c_2(E)$,
and ${{\mathcal E}xt}^1(E, {\mathcal O}_S) \cong {{\mathcal
E}xt}^2(\varkappa, {\mathcal O}_S).$ Then  ${{\mathcal F}itt}^0
{{\mathcal E}xt}^2(\varkappa, {\mathcal O}_S)$ is a sheaf of ideals
of (in general case nonreduced) subscheme  $Z$ of bounded length
\cite{Tim5} supported in a finite number of points on the surface
$S$. Hence, as it is shown in  \cite{Tim3}, the rest irreducible
components
 $\widetilde S_i, i>0$ of the scheme $\widetilde S$ in general case
can carry nonreduced scheme structure.

We assign to each semistable torsion-free coherent sheaf $E$ a pair
 $((\widetilde S, \widetilde L),
\widetilde E)$ with  $(\widetilde S, \widetilde L)$  defined as
described before.

Let $U$ be Zariski-open subset in one of components $\widetilde
S_i, i\ge 0$, and $\sigma^{\ast}E|_{\widetilde S_i}(U)$ be the
corresponding group of sections which is a  ${\mathcal
O}_{\widetilde S_i}(U)$-module. Sections $s\in
\sigma^{\ast}E|_{\widetilde S_i}(U)$ which are annihilated by
prime ideals of positive codimension in ${\mathcal O}_{\widetilde
S_i}(U)$, form a submodule in $\sigma^{\ast}E|_{\widetilde
S_i}(U)$. This submodule will be denoted as $tor\!s_i(U)$. The
correspondence $U \mapsto tor\!s_i(U)$ defines a subsheaf
$tor\!s_i \subset \sigma^{\ast}E|_{\widetilde S_i}.$ Note that
associated primes of positive codimension which annihilate
sections $s\in \sigma^{\ast}E|_{\widetilde S_i}(U)$, correspond to
subschemes supported in the preimage $\sigma^{-1}({\rm Supp\,}
\varkappa)=\bigcup_{i>0}\widetilde S_i.$ Since by the construction
the scheme  $\widetilde S=\bigcup_{i\ge 0}\widetilde S_i$ is
connected, then subsheaves  $tor\!s_i, i\ge 0,$ allow to form a
subsheaf  $tor\!s \subset \sigma^{\ast}E$. This subsheaf is
defined in the following way. A section  $s\in
\sigma^{\ast}E|_{\widetilde S_i}(U)$ satisfies the condition
$s\in tor\!s|_{\widetilde S_i}(U)$ if and only if
\begin{itemize}
\item{there exists a section $y\in {\mathcal O}_{\widetilde S_i}(U)$ such that $ys=0$,}
\item{at least one of following two conditions is satisfied: either  $y\in {\mathfrak p}$,
where $\mathfrak p$ is prime ideal of positive codimension; or there
exist a Zariski-open subset  $V\subset \widetilde S$ and a section
$s' \in \sigma^{\ast}E (V)$ such that $V\supset U$, $s'|_U=s$, and
$s'|_{V\cap \widetilde S_0} \in tor\!s (\sigma^{\ast}E|_{\widetilde
S_0})(V\cap \widetilde S_0)$. The torsion subsheaf
$tor\!s(\sigma^{\ast}E|_{\widetilde S_0})$ in the recent expression
is understood in the usual sense.}
\end{itemize}

In our construction the subsheaf $tor\!s \subset \sigma^{\ast}E$
plays the role which is completely analogous to the role of the
 subsheaf of torsion for the case of reduced and irreducible
 base scheme. Since no ambiguity occur, the symbol $tor\!s$
 anywhere in this article
is understood as mentioned above. The subsheaf $tor\!s$ is called a
{\it torsion subsheaf}.

In  \cite{Tim4} it is proven that sheaves $\sigma^{\ast} E/tor\!s$
are locally free. The sheaf  $\widetilde E$ in the pair
$((\widetilde S, \widetilde L), \widetilde E)$ is defined by the
formula $\widetilde E = \sigma^{\ast}E/tor\!s$. In this case there
is an isomorphism $H^0(\widetilde S, \widetilde E \otimes
\widetilde L) \cong H^0(S,E\otimes L).$

In \cite{Tim4E} it is shown that the restriction of the sheaf
 $\widetilde E$ on each component $\widetilde S_i$, $i>0,$ is given by the
 condition of quasi-ideality (\ref{quasiideal}) where $q_0: {\mathcal
O}_S^{\oplus r}\twoheadrightarrow \varkappa$ is an epimorphism
defined by the exact triple  $0\to E \to E^{\vee \vee} \to
\varkappa \to 0$ and by  local freeness of the sheaf $E^{\vee
\vee}$.

Resolution of singularities of a semistable sheaf  $E$ can be
globalized in a flat family by the procedure developed in articles
 \cite{Tim1,Tim2, Tim4} in different settings. Let  $T$ be
reduced irreducible  quasi-projective scheme, ${\mathbb E}$ be a
 sheaf of ${\mathcal O}_{T\times
S}$-modules, ${\mathbb L}$ invertible  ${\mathcal O}_{T\times
S}$-sheaf very ample relatively to $T$ and such that ${\mathbb
L}|_{t\times S}=L$ and  $\chi({\mathbb E} \otimes {\mathbb
L}^m|_{t\times S})=rp(m)$ for all closed points  $t\in T$. Also
suppose that  $T$ contains nonempty open subset $T_0$ such that
 ${\mathbb E}|_{T_0\times S}$ is locally free ${\mathcal O}_{T_0
 \times S}$-module. Then the following data is defined:
\begin{itemize} \item{$\widetilde T$ integral normal scheme obtained
by some blowing up  $\phi: \widetilde T \to T$ of the scheme $T$,}
\item{$\pi: \widetilde \Sigma \to \widetilde T$ flat family of
admissible schemes with invertible ${\mathcal O}_{\widetilde
\Sigma}$-module  $\widetilde {\mathbb L}$ such that $\widetilde
{\mathbb L}|_{t\times S}$ is the distinguished polarization of the
scheme  $\pi^{-1}(t)$,}
\item{$\widetilde {\mathbb E}$ locally free
${\mathcal O}_{\widetilde \Sigma}$-module and $((\pi^{-1}(t),
\widetilde {\mathbb L}|_{\pi^{-1}(t)}), \widetilde {\mathbb
E}|_{\pi^{-1}(t)})$ is $S$-semistable admissible pair.}
\end{itemize}
In this case there is a blowup morphism $\Phi: \widetilde \Sigma
\to
 \widetilde T \times S$, and $(\Phi_{\ast} \widetilde {\mathbb E})^{\vee
\vee}=\phi^{\ast}{\mathbb E}$. This follows from the coincidence
of sheaves from the  left hand side  and from the right hand side,
on the open subset out of a subscheme of codimension 3. Both
sheaves are reflexive. The scheme  $\widetilde T \times S$ is
integral and normal.

The mechanism described is called in \cite{Tim4} as  {\it standard
resolution}.

In section  1 we recall definitions of reduced functors
$(\mathfrak f^{GM}/\sim)$ of moduli of coherent semistable
torsion-free sheaves ("the Gieseker -- Maruyama moduli functor")
and $(\mathfrak f/\sim)$ of moduli of admissible semistable pairs.
The rank $r$ and the polynomial $p(m)$ are fixed and equal for
both moduli functors.
\smallskip

In the present article we prove following results.
\begin{theorem}\label{thfunc} There is a morphism of reduced moduli
$\mathfrak t: (\mathfrak f^{GM}/\!\sim)\to (\mathfrak f/\!\sim),$
defined by the procedure of standard resolution.
\end{theorem}
\begin{theorem}\label{thsch} Main components of the reduced scheme
$\widetilde M$ are isomorphic to main components of the reduced
Gieseker -- Maruyama scheme.
\end{theorem}
\smallskip
These theorems are proven in sections 1 and 2 respectively.

\section{Morphism of moduli functors: proof of theorem \ref{thfunc}}

Following  \cite[ch. 2, sect. 2.2]{HL} we recall some definitions.
Let ${\mathcal C}$ be a category, ${\mathcal C}^o$ its dual,
 ${\mathcal C}'={{\mathcal F}unct}({\mathcal C}^o, Sets)$ the category
 of functors to the category of sets. By Yoneda's lemma the functor
 ${\mathcal C} \to
{\mathcal C}': F\mapsto (\underline F: X\mapsto {\rm Hom
\,}_{{\mathcal C}}(X, F))$ includes  ${\mathcal C}$ in ${\mathcal
C}'$ as a full subcategory.

\begin{definition}\cite[ch. 2, definition 2.2.1]{HL} The functor
 ${\mathfrak f} \in
{{\mathcal O}}b\, {\mathcal C}'$ is {\it corepresented by the
object} $F \in {{\mathcal O}}b \,{\mathcal C}$ if there is a
${\mathcal C}'$-morphism $\psi : {\mathfrak f} \to \underline F$
such that  $\psi': {\mathfrak f} \to \underline F'$ factors
through the unique morphism $\omega: \underline F \to \underline
F'$.
\end{definition}

\begin{definition} The scheme  $\widetilde M$ is {\it coarse moduli space}
for the functor  $\mathfrak f$ if $\mathfrak f$ is corepresented
by the scheme
 $\widetilde M.$
 \end{definition}

Let  $T$ be a scheme over the field  $k$. Consider families of
semistable pairs \begin{equation}{\mathfrak F}_T=\left\{
\begin{array}{c}\widetilde \pi: \widetilde \Sigma_T \to T,  \;
\widetilde {\mathbb L}_T\in Pic \,\widetilde \Sigma_T ,
\;\forall t\in T \;\widetilde L_t=\widetilde {\mathbb L}_T|_{\widetilde \pi^{-1}(t)}\mbox{\rm \;\; very\;ample};\\
(\widetilde \pi^{-1}(t),\widetilde L_t) \mbox{\rm \;admissible scheme with distinguished  }\\
\mbox{\rm polarization}; \\
 \widetilde {\mathbb E}_T - \mbox{\rm locally free  } {\mathcal O}_{\widetilde \Sigma_T}-\mbox{\rm
 sheaf};\\
 \chi(\widetilde {\mathbb E}_T\otimes\widetilde {\mathbb L}_T^{m})|_{\pi^{-1}(t)})=
 rp(m);\\
 ((\widetilde \pi^{-1}(t), \widetilde L_t), \widetilde {\mathbb E}_T|_{\widetilde \pi^{-1}(t)}) -
 \mbox{\rm (semi)stable admissible pair}
 \end{array} \right\} \nonumber\end{equation} and a functor
 $\mathfrak f: (Schemes_k)^o \to (Sets)$ from the category of
 $k$-schemes to the category of sets, assigning to a scheme $T$ the set
$\{\mathfrak F_T\}$. The moduli functor $(\mathfrak f/\sim)$
assigns to a scheme  $T$ the set of equivalence classes
$(\{\mathfrak F_T\}/\sim)$.

The equivalence relation $\sim$ is defined as follows. Families
 $((\widetilde \pi: \widetilde \Sigma_T \to T, \widetilde {\mathbb L}_T), \widetilde {\mathbb E}_T)$
 and
 $((\widetilde \pi': \widetilde \Sigma'_T \to T, \widetilde {\mathbb L}'_T), \widetilde
 {\mathbb E}'_T)$ from the set $\{{\mathfrak F}_T\}$ are said to be equivalent  (notation:
 $((\widetilde \pi: \widetilde \Sigma_T \to T,
 \widetilde {\mathbb L}_T), \widetilde {\mathbb E}_T) \sim
 ((\widetilde \pi': \widetilde \Sigma'_T \to T, \widetilde {\mathbb L}'_T), \widetilde
 {\mathbb E}'_T)$) if\\
 1) there is an isomorphism $\widetilde \Sigma_T
 \stackrel{\sim}{\longrightarrow} \widetilde \Sigma'_T$ such that the diagram
 \begin{equation*}\xymatrix{\widetilde \Sigma_T \ar[rd]_{\pi}\ar[rr]^{\sim}&&\widetilde \Sigma'_T \ar[ld]^{\pi'}\\
&T }
 \end{equation*}
 commutes;\\
 2) there exist  linear bundles $\mathcal L', \mathcal L''$ on
 $T$ such that  $\widetilde {\mathbb E}'_T = \widetilde
 {\mathbb E}_T
 \otimes \widetilde \pi^{\ast} \mathcal
 L'$ and $\widetilde
 \L'_T=\widetilde \L_T \otimes \widetilde \pi'^{\ast} \mathcal L''.$

\begin{convention} \label{redconv} In this paper we restrict ourselves
by the consideration of the full subcategory  $RSch_k$ of reduced
schemes and of the  {\it reduced moduli functor} $(\mathfrak
f_{red}/\sim)=(\mathfrak f|_{(RSch_k)^o}/\sim ) $ \cite{Tim5}.
Since no ambiguity occur, we use the symbol $(\mathfrak f/\sim)$
for the reduced moduli functor.
\end{convention}

In the general case results of the paper \cite{Tim4} provide
existence of a coarse moduli space for any maximal under inclusion
irreducible substack in  $\coprod_{T\in {\mathcal O} b(RSch_k)}
({\mathfrak F}_T/\!\sim)$ if this substack contains such pairs
$((\widetilde \pi^{-1}(t), \widetilde L_t), \widetilde {\mathbb
E}_T|_{\pi^{-1}(t)})$ that $(\pi^{-1}(t), \widetilde L_t)\cong
(S,L)$. Such pairs will be referred to as {\it $S$-pairs}. We
mention under  $\widetilde M$ namely the moduli space of a
substack containing semistable $S$-pairs and emphasize this
speaking about {\it main components} of the moduli scheme.

Analogously, we mention the Gieseker -- Maruyama scheme $\overline
M$ as union of those components of reduced moduli scheme of
semistable torsion-free sheaves, that contain locally free
sheaves.

\smallskip

The Gieseker -- Maruyama functor $\mathfrak f^{GM}: (Schemes_k)^o
\to Sets$ is defined as follows: $T\mapsto \{ \mathfrak
F_T^{GM}\}$, where \begin{equation*} \mathfrak F_T^{GM}= \left\{
\begin{array}{c} \E_T \mbox{\rm  is a sheaf of } \OO_{T\times S}-
\mbox{\rm modules, flat over } T;\\
\L_T \mbox{\rm  is invertible sheaf  of } \OO_{T\times
S}-\mbox{\rm  modules,}\\
L_t:=\L_T|_{t\times S} \;\mbox{\rm is very ample;}
\\
E_t:=\E_T|_{t\times S} \;\mbox{\rm is torsion-free }\\
\mbox{ \rm and Gieseker-semistable with respect to  } L_t;\\
\chi(E_t \otimes L_t^m)=rp(m).\end{array}\right\}
\end{equation*}

Families  $(\L_T,\E_T)$ and  $(\L'_T,\E'_T)$ are said to be
equivalent if there are invertible $\OO_T$-sheaves $\LL'$ and
$\LL''$ such that for the projection $p: T\times S\to T$ one has
\begin{eqnarray*}
\E'_T&\cong &\E_T \otimes p^{\ast}\LL',\\
\L'_T&\cong &\L_T \otimes p^{\ast}\LL''.
\end{eqnarray*}

For this functor we use  the convention which is analogous to the
convention \ref{redconv}.

The functor morphism  $\mathfrak t: (\mathfrak f^{GM}/\sim) \to
(\mathfrak f/\sim )$ is defined by commutative diagrams
\begin{equation}\label{morfun}\xymatrix{T \ar@{|->}[rd] \ar@{|->}[r]&
\{\mathfrak F^{GM}_T\}/\sim \ar[d]^{\mathfrak t_T}\\
& \{ \mathfrak F_T\}/\sim}
\end{equation}
where  $T\in Ob RSch_k$, $\mathfrak t_T:(\{\mathfrak
F^{GM}_T\}/\sim) \to (\{ \mathfrak F_T\}/\sim)$ is a morphism in
the category of sets (the map of sets).

\begin{remark}\label{deformability} We consider subfunctors in $\mathfrak f^{GM}$ and in
$\mathfrak f$ which correspond to maximal under inclusion
irreducible substacks containing locally free sheaves and
$S$-pairs respectively. Hence any family  $\mathfrak F^{GM}_{T}$
(respectively, $\mathfrak F_T$) can be include in a family
$\mathfrak F^{GM}_{T'}$ (respectively, $\mathfrak F_{T'}$) with
some connected base $T'$ and containing locally free sheaves
(respectively, $S$-pairs) according to the fibres diagram
\begin{equation*}\label{deffam}\xymatrix{T \ar@{^(->}[d]^i & \ar[l] \mathfrak F^{GM}_T \ar@{^(->}[d] \\
T' & \ar[l] \mathfrak F^{GM}_{T'}} (\mbox{\rm respectively, }
\xymatrix{T \ar@{^(->}[d]^i & \ar[l] \mathfrak F_T \ar@{^(->}[d] \\
T' & \ar[l] \mathfrak F_{T'}})
\end{equation*}
Namely, $\E_T=i^{\ast}\E_{T'}$ (respectively, $\widetilde
\Sigma_T= \widetilde \Sigma_{T'} \times _{T'} T$, $\widetilde i:
\widetilde \Sigma_T \hookrightarrow \widetilde \Sigma_{T'}$ is the
induced morphism of immersion, $\widetilde \E_T = \widetilde
i^{\ast} \widetilde \E_{T'},$ $\widetilde \L_T= \widetilde
i^{\ast} \widetilde \L_{T'}$). In particular, this restriction
excludes embedded components of moduli scheme from our
consideration whenever these components do not contain locally
free sheaves (respectively, $S$-pairs). Then it is enough to
construct diagrams (\ref{morfun}) only for families which contain
locally free sheaves (respectively, $S$-pairs), where $T$ is
reduced scheme.\end{remark}

Let  $p: \Sigma_T \to T$ be a flat family of schemes such that its
fibre is isomorphic to $S$, $\L_T$ be a family of very ample
invertible sheaves on fibres of the family $p$, $\E_T$ be a flat
family of coherent torsion-free sheaves on fibres of $p$. The
sheaves are mentioned to have rank  $r$ and Hilbert polynomial
$rp(m)$ and to be semistable with respect to polarizations induced
by the family  $\L_T$. The application of the standard resolution
leads to the collection of data $(\widetilde \pi: \widetilde
\Sigma_{\widetilde T} \to \widetilde T, \widetilde \L_{\widetilde
T}, \widetilde \E_{\widetilde T}).$ Let  $\Sigma_{\widetilde
T}:=\Sigma_T \times _T \widetilde T$ where  $\phi: \widetilde T
\to T$ is birational morphism also provided by the standard
resolution, and $ (\phi,id_S): \Sigma_{\widetilde T} \to \Sigma_T$
is the induced morphism.

Further, due to the considerations of \cite{Tim1,Tim2,Tim4} there
is a  {\it partial} morphism of functors \linebreak $\mathfrak
t^{-1}: (\mathfrak f/\sim) \to (\mathfrak f^{GM}/\sim)$ defined by
the morphism  $\sigma \!\!\! \sigma: \widetilde \Sigma_{\widetilde
T} \to \Sigma_{\widetilde T}$ and the operation $(\sigma \!\!\!
\sigma_{\ast} -)^{\vee \vee}$ on those families which are obtained
by standard resolution from families of coherent semistable
torsion-free sheaves. Then  $\mathfrak t^{-1} \circ \mathfrak
t=id_{Sets}.$ Since  $\mathfrak t^{-1}$ is defined only partially,
it is impossible to claim that  $\mathfrak t$ is an isomorphism.

\begin{remark} Also, as it is shown in \cite{Tim4}, there is
a birational morphism of moduli schemes $\kappa: \overline M \to
\widetilde M$ which are mentioned to be reduced \cite{Tim4E}. The
scheme  $\widetilde M$ can be not normal. Hence, although $\kappa$
is a bijective morphism and becomes a morphism of integral schemes
if restricted on each of irreducible components, this does not
imply that  $\kappa$ is an isomorphism.
\end{remark}

In the further text we will show that there is a morphism of
reduced Gieseker -- Maruyama moduli functor to the reduced moduli
functor of admissible semistable pairs. Namely, for any reduced
scheme $T$ the correspondence $\E_T \mapsto ((\widetilde \Sigma_T,
\widetilde \L_T), \widetilde \E_T)$ will be constructed. This
correspondence defines the map of sets $(\{\E_T\}/\sim) \to
(\{((\widetilde \Sigma_T, \widetilde \L_T), \widetilde
\E_T)\}/\sim)$. This all means that for any family of semistable
coherent sheaves $\E_T$ which is flat over its base $T$ one can
build up a family $((\widetilde \Sigma_T, \widetilde \L_T),
\widetilde \E_T)$ of admissible semistable pairs with the same
base $T$.

The procedure of standard resolution yields in the family of
admissible schemes $\widetilde \pi:\widetilde \Sigma_{\widetilde
T} \to \widetilde T$ which is flat over $\widetilde T$, in the
locally free $\OO_{\widetilde \Sigma_{\widetilde T}}$-sheaf
$\widetilde \E_{\widetilde T}$ and in the invertible
$\OO_{\widetilde \Sigma_{\widetilde T}}$-sheaf $\widetilde
\L_{\widetilde T}$, which is very ample with respect to the
morphism  $\widetilde \pi$.

\begin{proposition} \label{descent} There exist \begin{itemize}
\item{$\widetilde \Sigma_T$ scheme,}
\item{$\pi: \widetilde \Sigma_T \to T$ flat morphism,}
\item{$\overline \phi: \widetilde \Sigma_{\widetilde T} \to
\widetilde \Sigma_T$ birational morphism,}
\item{$\widetilde \E_T$ locally free $\OO_{\widetilde \Sigma_T}$-sheaf,}
\item{$\widetilde \L_T$ invertible $\OO_{\widetilde \Sigma_T}$-sheaf,}
\end{itemize}
such that  \begin{itemize}\item{the square
\begin{equation*}\xymatrix{\widetilde \Sigma_{\widetilde T}
\ar[d]_{\widetilde \pi} \ar[r]^{\overline \phi}&
\widetilde \Sigma_T \ar[d]^{\pi}\\
\widetilde T \ar[r]^{\phi}& T}
\end{equation*} is Cartesian,}
\item{$\widetilde \E_{\widetilde T} \otimes \LL'=
\overline \phi^{\ast}\widetilde \E_T$ for some $\LL' \in \Pic
\widetilde T$,}
\item{$\widetilde \L_{\widetilde T} \otimes \LL''=
\overline \phi^{\ast} \widetilde \L_T$ for some $\LL'' \in \Pic
\widetilde T.$}
\end{itemize}
\end{proposition}
The proposition formulated implies the functor morphism of
interest $\mathfrak t:(\mathfrak f^{GM}/\sim) \to (\mathfrak
f/\sim)$. It is defined for any reduced scheme  $T \in Ob\,
RSch_k$ by the commutative diagram (\ref{morfun})
\begin{equation*} \xymatrix{T \ar@{|->}[rd]\ar@{|->}[r]
& (\{\E_T\}/\sim)\ar[d]\\
& (\{((\widetilde \Sigma_T, \widetilde \L_T), \widetilde
\E_T)\}/\sim)}
\end{equation*}
where the right vertical arrow is a morphism (mapping) in the
category of sets. This mapping is defined by the proposition
\ref{descent}. The horizontal and the skew arrows are defined by
functorial correspondences $(\mathfrak f^{GM}/\sim)$ and
$(\mathfrak f/\sim)$ respectively.

\begin{proof}[Proof of proposition \ref{descent}.] For the construction
of the scheme  $\widetilde \Sigma_T$ we assume that $m\gg 0$ is as
big as the sheaf of $\OO_{\widetilde T}$-modules $\widetilde
\pi_{\ast}(\widetilde \E_{\widetilde T} \otimes \widetilde
\L_{\widetilde T}^m)$ is locally free, the canonically defined
morphism $\widetilde \pi^{\ast} \widetilde \pi_{\ast} (\widetilde
\E_{\widetilde T} \otimes \widetilde \L_{\widetilde T}^m )\to
\widetilde \E_{\widetilde T} \otimes \widetilde \L_{\widetilde
T}^m$ is surjective, and there is in induced closed immersion
$\widetilde \Sigma_{\widetilde T} \hookrightarrow G(\widetilde
\pi_{\ast} (\widetilde \E_{\widetilde T} \otimes \widetilde
\L_{\widetilde T}^m), r)$ into Grassmannian bundle $G(\widetilde
\pi_{\ast} (\widetilde \E_{\widetilde T} \otimes \widetilde
\L_{\widetilde T}^m), r)$. Also there is a (relative) Pl\"{u}cker
immersion of the Grassmannian bundle into the (relative)
projective space
$$G(\widetilde \pi_{\ast} (\widetilde \E_{\widetilde T} \otimes
\widetilde \L_{\widetilde T}^m), r) \hookrightarrow P(\bigwedge^r
\widetilde \pi_{\ast} (\widetilde \E_{\widetilde T} \otimes
\widetilde \L_{\widetilde T}^m)).$$

Besides consider the  isomorphism of $\OO_{\widetilde T}$-sheaves
$p_{\ast} (\phi, id_S)^{\ast} (\E_T \otimes \L^m_T)= \phi^{\ast}
p_{\ast} (\E_T \otimes \L^m_T)$ and the sheaf  $\widetilde
\pi_{\ast} (\widetilde \E_{\widetilde T} \otimes \widetilde
\L_{\widetilde T}^m)$. These sheaves are locally free and coincide
on open subset of the scheme  $\widetilde T$. Then, if  $\LL'$
denotes the invertible sheaf of the form $\LL':= \det p_{\ast}
(\phi, id_S)^{\ast} (\E_T \otimes \L^m_T) \otimes (\det \widetilde
\pi_{\ast} (\widetilde \E_{\widetilde T} \otimes \widetilde
\L_{\widetilde T}^m))^{\vee}$, then there are two locally free
 $\OO_{\widetilde T}$-sheaves
$\widetilde \pi_{\ast} (\widetilde \E_{\widetilde T} \otimes
\widetilde \L_{\widetilde T}^m) \otimes \LL'$ and $p_{\ast} (\phi,
id_S)^{\ast} (\E_T \otimes \L^m_T).$ They coincide along the open
subset of the normal integral scheme $\widetilde T$. This subset
is obtained by excision of some subscheme of codimension $\ge 2$.
This implies that these sheaves coincide, namely,
\begin{equation*} \widetilde \pi_{\ast} (\widetilde \E_{\widetilde T}
\otimes \widetilde \L_{\widetilde T}^m) \otimes \LL'= p_{\ast}
(\phi, id_S)^{\ast} (\E_T \otimes \L^m_T)= \phi^{\ast} p_{\ast}
(\E_T \otimes \L^m_T).
\end{equation*}
Consequently, formation of exterior powers and passing to
projectivizations and to Grassmannian bundles induces the fibred
diagram
\begin{equation}\label{cart1}
\xymatrix{P(\bigwedge^r (\widetilde \pi_{\ast} (\widetilde
\E_{\widetilde T} \otimes \widetilde \L_{\widetilde T}^m) \otimes
\LL')) \ar[d] \ar[r]& P(\bigwedge^r p_{\ast} (\E_T
\otimes \L^m_T)) \ar[d]\\
\widetilde T \ar[r]& T}
\end{equation}
It is compatible with Pl\"{u}cker embeddings in the sense that the
square \begin{equation}\label{cart2}\xymatrix{P(\bigwedge^r
(\widetilde \pi_{\ast} (\widetilde \E_{\widetilde T} \otimes
\widetilde \L_{\widetilde T}^m) \otimes \LL')) \ar[r]&
P(\bigwedge^r p_{\ast} (\E_T
\otimes \L^m_T)) \\
G(\widetilde \pi_{\ast} (\widetilde \E_{\widetilde T} \otimes
\widetilde \L_{\widetilde T}^m) \otimes \LL'),r)\ar@{^(->}[u]
\ar[r]^g &G(p_{\ast} (\E_T \otimes \L^m_T),r) \ar@{^(->}[u]}
\end{equation}
is also Cartesian. Denote by $\widetilde \Sigma'_T$ the
(scheme-theoretic) image of the composite map $\widetilde
\Sigma_{\widetilde T} \hookrightarrow G(\widetilde \pi_{\ast}
(\widetilde \E_{\widetilde T} \otimes \widetilde \L_{\widetilde
T}^m) \otimes \LL'),r) \stackrel{g}{\longrightarrow}G(p_{\ast}
(\E_T \otimes \L^m_T),r).$ The immersion here is induced by the
sheaf  $\widetilde \pi_{\ast} (\widetilde \E_{\widetilde T}
\otimes \widetilde \L_{\widetilde T}^m) \otimes \LL'$; the
morphism  $g$ comes from the previous diagram. To convince that
the scheme  $\widetilde \Sigma'_T$ is  flat over $T$, it is enough
to verify the coincidence of images for fibres of the scheme
$\widetilde \Sigma_{\widetilde T}$ over closed points of the fibre
$\phi^{-1}(t)$ for each closed point $t\in T$. But by the
construction, images of all fibres of the scheme $\widetilde
\Sigma_{\widetilde T}$ over points of the subscheme
$\phi^{-1}(t)$, have equal collections of local equations in
fibres of the projective bundle $P(\bigwedge^r p_{\ast} (\E_T
\otimes \L^m_T))$. Consequently, they  are also given by equal
local equations in fibres of the Grassmannian bundle $G(p_{\ast}
(\E_T \otimes \L^m_T),r)$. This completes the proof of flatness of
 $\widetilde
\Sigma'_{T}$ over $T$.

Also it is clear that the morphism $\overline \phi': \widetilde
\Sigma_{\widetilde T} \to \widetilde \Sigma'_T$ of the scheme
$\widetilde \Sigma_{\widetilde T}$ onto its image becomes an
isomorphism when restricted to the open subset $\pi^{-1}(T_0)$.
$T_0$ is open subscheme of the scheme $T$ where $\phi$ is an
isomorphism. This implies that the morphism $\overline \phi'$ is
birational.

Now turn to closed immersions $j_{\widetilde T}: \widetilde
\Sigma_{\widetilde T} \hookrightarrow G(\widetilde \pi_{\ast}
(\widetilde \E_{\widetilde T} \otimes \widetilde \L_{\widetilde
T}^m) \otimes \LL'),r)$ and $j:\widetilde \Sigma'_T
\hookrightarrow G(p_{\ast} (\E_T \otimes \L^m_T),r).$  We denote
by the symbol $\SS_T$ the universal (locally free) quotient sheaf
of rank $r$ on $G(p_{\ast} (\E_T \otimes \L^m_T),r)$. Also the
symbol $\SS_{\widetilde T}$ denotes the universal (locally free)
quotient sheaf of rank $r$  on $G(\widetilde \pi_{\ast}
(\widetilde \E_{\widetilde T} \otimes \widetilde \L_{\widetilde
T}^m) \otimes \LL'),r)$. It is clear that by the construction
(\ref{cart1}, \ref{cart2}) of Grassmannian bundles of our interest
we can write $\SS_{\widetilde T}=g^{\ast} \SS_T$. Then
$j_{\widetilde T}^{\ast} \SS_{\widetilde T}=\widetilde
\E_{\widetilde T} \otimes \widetilde \L_{\widetilde T}^{m} \otimes
\widetilde \pi^{\ast} \LL'$, and $j_T^{\ast} \SS_T$ is a locally
free sheaf on the scheme $\widetilde \Sigma'_{T}$.

Consider the invertible  $\OO_{\widetilde \Sigma_{\widetilde
T}}$-sheaf $\widetilde \L_{\widetilde T}$ providing the
distinguished polarization on fibres of the morphism $\widetilde
\pi: \widetilde \Sigma_{\widetilde T} \to \widetilde T$, and its
direct image  $ \widetilde \pi_{\ast} \widetilde \L_{\widetilde
T}$. The recent sheaf is locally free by the choice of the sheaf
 $\widetilde
\L_{\widetilde T}$. Also take (any) invertible
$\OO_{\Sigma_T}$-sheaf $\L_T$ which is very ample relatively to
the projection $p: \Sigma_{T} \to T$ and such that sheaves
$\widetilde \L_{\widetilde T}$  and $\L_{T}$ induce equal
polarizations on fibres over points of open subscheme $T_0$. The
sheaf $p_{\ast} \L_T$ is also locally free. Since the Hilbert
polynomial on fibres of morphisms $\widetilde \pi$ and $p$ is
constant over the base and the morphism  $\widetilde
\Sigma'_{\widetilde T} \to \Sigma_T$ is birational, ranks of
locally free sheaves  $\widetilde \pi_{\ast} \widetilde
\L_{\widetilde T}$ and $p_{\ast} \L_T$ are equal. Moreover, their
restrictions on the open subscheme  $T_0 \subset T$  are
isomorphic. Then there exists an invertible $\OO_{\widetilde
T}$-sheaf $\LL''$ such that  $\widetilde \pi_{\ast} \widetilde
\L_{\widetilde T} \otimes \LL ''=\phi^{\ast} p_{\ast}
\L_T=p_{\ast} (\phi, id_S)^{\ast} \L_T$.

Now consider the relative projective space $P(\widetilde
\pi_{\ast} \widetilde \L_{\widetilde T}\otimes \LL '')$ together
with closed immersion $i: \widetilde \Sigma_{\widetilde
T}\hookrightarrow P(\widetilde \pi_{\ast} \widetilde
\L_{\widetilde T}\otimes \LL '')$. Also take the relative
projective space of same dimension $P(p_{\ast} \L_T)$. Then the
square
\begin{equation*} \xymatrix{
P(\widetilde \pi_{\ast} \widetilde \L_{\widetilde T} \otimes \LL'')
\ar[d] \ar[r]& P(p_{\ast}\L_T) \ar[d]\\
\widetilde T \ar[r]& T}
\end{equation*}
is Cartesian. Denote by  $\widetilde \Sigma''_T$ the image of the
composite map $ \widetilde \Sigma_{\widetilde T} \hookrightarrow
P(\widetilde \pi_{\ast} \widetilde \L_{\widetilde T} \otimes
\LL'') \to P(p_{\ast} \L_T)$. By the construction the scheme
$\widetilde \Sigma''_{T}$ is flat over $T$. Denote by $j:
\widetilde \Sigma ''_T \hookrightarrow P(p_{\ast} \L_T)$ the
corresponding closed immersion of this subscheme; let $\OO_P(1)$
be the canonical invertible sheaf on the projective bundle
$P(p_{\ast}\L_T)$. Then define $\widetilde \L''_T:=j^{\ast}
\OO_P(1)$.

Now form the fibred product $G(p_{\ast} (\E_T \otimes
L_T^m),r)\otimes _T P(p_{\ast}\L_T)$ and consider the mapping
$\widetilde \Sigma_{\widetilde T} \to G(p_{\ast} (\E_T \otimes
L_T^m),r)\times _T P(p_{\ast}\L_T)$. This map is induced by
mappings into each factor constructed earlier. Let  $\widetilde
\Sigma_T$ be the scheme image of this map. Then there is the
following commutative diagram:
\begin{equation*}\xymatrix{
\widetilde \Sigma'_T \ar@{^(->}[d] &\ar[l]_{p'} \widetilde \Sigma_T
\ar@{^(->}[d]
\ar[r]^{p''} & \widetilde \Sigma''_T \ar@{^(->}[d]\\
G(p_{\ast} (\E_T \otimes \L_T^m),r) & \ar[l] G(p_{\ast} (\E_T
\otimes \L_T^m),r) \times _T P(p_{\ast} \L_T)  \ar[r]&
P(p_{\ast}\L_T)}
\end{equation*}
Birational morphisms  $p', p''$ are defined by projections of the
product  $G(p_{\ast} (\E_T \otimes \L_T^m),r) \times _T P(p_{\ast}
\L_T)$ to each factor.

It rests to define sheaves  $\widetilde \L_T$ and $\widetilde
\E_T$ on the scheme  $\widetilde \Sigma_T$ by formulas:
\begin{eqnarray*}
\widetilde \L_T&:=&p''^{\ast} \widetilde \L''_T,\\
\widetilde \E_T&:=&p'^{\ast} (j_T^{\ast} \SS_T) \otimes \widetilde
L_T^{-m}
\end{eqnarray*}
  \end{proof}
\begin{remark} If $T$ is an integral and normal scheme then
the functorial correspondence which has been constructed is
invertible. Then the class of schemes where the partial functor
$\mathfrak t^{-1}$ is defined, contains all integral and normal
schemes.
\end{remark}

\section{Isomorphism of moduli schemes: proof of theorem \ref{thsch}}

First recall some objects and constructions which are used in
classical built-up of Gieseker -- Maruyama scheme. By choice of
the polarization  $L$ we assume that $H^0(S, E \otimes L^m)
\otimes L^{\vee} \to E$ is an epimorphism. Let ${\rm
Quot\,}^{rp(m)}(V \otimes L^{-m})$ be the Grothendieck's $Quot$
scheme. It parameterizes quotient sheaves $V \otimes
L^{-m}\twoheadrightarrow E$ for $V \cong H^0(S, E \otimes L^m),$
with Hilbert polynomial $rp(m)$ with respect to $L$. The scheme
$Quot$ carries universal family of quotient sheaves
\begin{equation*}
V \otimes L^{-m} \boxtimes {\mathcal O}_{{\rm Quot\,}^{rp(m)}(V
\otimes L^{-m})}\twoheadrightarrow {\mathbb E}_{{\rm Quot\,}}.
\end{equation*}
Let $Q'$ be the quasiprojective subscheme in ${\rm
Quot\,}^{rp(m)}(V \otimes L^{-m}),$ consisting of points
corresponding to semistable sheaves, $\xi: Q\to Q'$ be any its
smooth resolution, ${\mathbb E}_Q:=\xi^{\ast}{\mathbb E}_{{\rm
Quot\,}}|_{Q'}$ family of coherent sheaves which is flat over $Q$.
It comes from ${\mathbb E}_{{\rm Quot\,}}$. Classical way to
construct Gieseker -- Maruyama scheme is to form a GIT-quotient
$Q'^{ss}// PGL(V)$ of the set  $Q'^{ss}$ of semistable points with
respect to the action of the group $PGL(V)$ upon the scheme $Q'$.
The action is induced by linear transformations of  vector space
$V$.

Application of standard resolution to the data $Q\times S,
{\mathbb L}={\mathcal O}_Q \boxtimes L, {\mathbb E}_Q$ leads to
the collection of data $\widetilde Q, \widetilde {\mathbb L},
\widetilde {\mathbb E}_Q$.

Note that the sheaf  $\widetilde E \otimes \widetilde L^m$ defines
a closed immersion $j: \widetilde S \hookrightarrow G(V,r)$ where
 $G(V,r)$ is Grassmann variety parameterizing quotient spaces of
 dimension $r$ of the  vector space $V\cong H^0(\widetilde S,
\widetilde E \otimes \widetilde L^m)$. The immersion $j$ is
defined non-uniquely but up to isomorphy class $H^0(\widetilde S,
\widetilde E \otimes \widetilde L^m) \stackrel{\sim}
{\longrightarrow} V$ modulo multiplication by nonzero scalars
$\vartheta \in k^{\ast}.$ Let  $P(m)=\chi(j^{\ast}{\mathcal
O}_{G(V,r)}(m))$ be the Hilbert polynomial of the subscheme
$j(\widetilde S) \subset G(V,r)$. Hence the point corresponding to
the subscheme $j(\widetilde S) \subset G(V,r),$ is defined in the
Hilbert scheme  ${\rm Hilb \,}^{P(t)}G(V,r)$ up to the action of
the group  $PGL(V).$

In \cite{Tim4} it is shown that the data $\widetilde Q, \widetilde
\L, \widetilde {\mathbb E}_Q$ defines the morphism $\mu:\widetilde
Q \to {\rm Hilb \,}^{P(m)}G(V, r)$. Also it is proven there that
there is a morphism of projective schemes $\kappa: \overline M \to
\widetilde M$.

\smallskip
For the further consideration we also need the notion of
M-equivalence for semistable admissible pairs. This notion is
introduced and investigated in  \cite{Tim4}.

For any two schemes of the form $\widetilde S_1={\rm Proj \,}
\bigoplus_{s\ge 0}(I_1[t]+(t))^{s}/(t)^{s+1}$ и $\widetilde
S_2={\rm Proj \,} \bigoplus_{s\ge 0}(I_2[t]+(t))^{s}/(t)^{s+1}$
define a scheme
 $$\widetilde S_{12}={\rm Proj \,} \bigoplus_{s\ge
0}(I'_1[t]+(t))^{s}/(t)^{s+1}={\rm Proj \,} \bigoplus_{s\ge
0}(I'_2[t]+(t))^{s}/(t)^{s+1}$$ together with morphisms
$\widetilde S_1\stackrel{\sigma'_1}{\longleftarrow} \widetilde
S_{12} \stackrel{\sigma'_2}{\longrightarrow} \widetilde S_2,$ such
that the diagram \begin{equation*}\xymatrix{\widetilde S_{12}
\ar[r]^{\sigma'_2}\ar[d]_{\sigma'_1}& \widetilde S_2
\ar[d]^{\sigma_2}\\
\widetilde S_1 \ar[r]_{\sigma_1}&S}
\end{equation*}
commutes. The operation
 $(\widetilde S_1, \widetilde S_2)
\mapsto \widetilde S_1 \diamond \widetilde S_2=\widetilde S_{12}$
defined as before, is associative. Moreover, since for any
admissible morphism  $\sigma :\widetilde S\to S$ the relation
$\widetilde S \diamond S =S \diamond \widetilde S =\widetilde S$
holds, admissible morphisms for each class $[E]$ of S-equivalent
semistable coherent sheaves, generate a commutative monoid
$\diamondsuit[E]$ with binary operation $\diamond$ and neutral
element  $id_S: S \to S.$

\begin{definition}\cite{Tim4} Semistable pairs  $(\widetilde S, \widetilde E)$
and $(\widetilde S', \widetilde E')$ are {\it $M$-equivalent
(monoidally equivalent)} if for morphisms of $\diamond$-product
$\widetilde S \diamond \widetilde S'$ to both factors $\overline
\sigma': \widetilde S \diamond \widetilde S' \to \widetilde S$ and
$\overline \sigma: \widetilde S \diamond \widetilde S' \to
\widetilde S'$ and for associated polystable sheaves $\bigoplus_i
gr_i (\widetilde E)$ and $\bigoplus_i gr_i (\widetilde E')$ there
are isomorphisms
\begin{equation*}\overline \sigma'^{\ast}
\bigoplus_i gr_i (\widetilde E)/tor\!s \cong \overline \sigma^{\ast}
\bigoplus_i gr_i (\widetilde E')/tor\!s.\end{equation*}
\end{definition}

It follows from results of \cite{Tim4} that classes of
M-equivalence of admissible semistable pairs $((\widetilde S,
\widetilde L), \widetilde E)$ are bijective to classes of
S-equivalence of semistable coherent sheaves. The bijection is
established by the correspondence $ E \mapsto ((\widetilde S,
\widetilde L),\sigma^{\ast} \widetilde E/tor\!s)$. Then it is
enough to convince ourselves that the Gieseker -- Maruyama moduli
scheme
 $\overline M$ corepresents the functor $\mathfrak f$ as well
 as the scheme  $\widetilde M.$ Choose an object
 $((\widetilde S, \widetilde L), \widetilde E)$. By bijectivity
 of the set of classes of M-equivalence for semistable admissible
 pairs to the set of classes of S-equivalence of semistable
 coherent sheaves, the object  $((\widetilde S, \widetilde L), \widetilde E)$
defines a morphism  $h\in {\rm Hom \,} ({\rm Spec \,} k, \overline
M).$

Inversely, the morphism  $h\in {\rm Hom \,} ({\rm Spec \,} k,
\overline M)$ distinguishes a point which represents the class of
M-equivalent objects  $((\widetilde S, \widetilde L), \widetilde
E)$.

We construct  for any reduced scheme  $B$ and for a natural
transformation $\psi': \mathfrak f \to \underline F'$ the unique
natural transformation $\omega: \underline{\overline M}\to
\underline F'$ such that $\psi'=\omega \circ \psi.$ The natural
transformation  $\psi'$ corresponds to a flat family
 $\pi: \widetilde \Sigma \to B$ with fibrewise polarization
 $\widetilde {\mathbb L}$ and equipped with the family of locally
 free sheaves  $\widetilde {\mathbb E}$. It is enough to restrict
 ourselves by the case when $B$ is integral scheme. Indeed, any
 reduced scheme can be represented as union of irreducible
 components. Also we take into account the remark
 \ref{deformability}

Restriction on any fibre $\pi^{-1}(b)$ of the morphism  $\pi$
yields in an object $((\pi^{-1}(b), \widetilde {\mathbb
L}|_{\pi^{-1}(b)}), \widetilde {\mathbb E}|_{\pi^{-1}(b)})$ of the
class  $\mathfrak F$. Then there is a morphism into the
Grassmannian bundle $\widetilde \Sigma  \to G(\pi_{\ast}
(\widetilde {\mathbb E} \otimes \widetilde {\mathbb L}^m) , r)$.
This morphism becomes an immersion $j_b: \pi^{-1}(b)
\hookrightarrow G(H^0(\pi^{-1}(b), \widetilde {\mathbb E} \otimes
\widetilde {\mathbb L}^m|_{\pi^{-1}(b)}),r)$ when restricted to
fibres of the morphism $\pi$.

The sheaf  $\pi_{\ast} (\widetilde {\mathbb E} \otimes \widetilde
{\mathbb L}^m)$ is locally free and then the Grassmannian bundle
 $G(\pi_{\ast}
(\widetilde {\mathbb E} \otimes \widetilde {\mathbb L}^m) , r)$ is
locally trivial over $B$. Let  $\bigcup_i B_i = B$ be Zariski-open
trivializing cover. Subfamilies
 $\widetilde
\Sigma_i$ are given by fibred squares
\begin{equation*} \xymatrix{\widetilde \Sigma \ar[d]_{\pi} & \ar@{_(->}[l] \widetilde \Sigma_i \ar[d]^{\pi_i}\\
B & \ar@{_(->}[l] B_i}
\end{equation*}
where horizontal arrows are open immersions. Fix isomorphisms of
trivialization $\tau_i:G(\pi_{\ast} (\widetilde {\mathbb E}
\otimes \widetilde {\mathbb L}^m) , r)|_{B_i}\to G(V, r) \times
B_i.$ The composite map $\widetilde \Sigma_i \stackrel{\!j_i}{\to}
G(\pi_{\ast} (\widetilde {\mathbb E} \otimes \widetilde {\mathbb
L}^m) , r)|_{B_i}\stackrel{\!\!\tau_i}{\to} G(V, r) \times B_i
\stackrel{pr_1}{\to} G(V,r)$ defines a morphism of the base into
Hilbert scheme  $\mu_i: B_i \to {\rm Hilb \,}^{P(t)}G(V,r)$. The
morphism $\mu_i$ factors through the subscheme $\mu(\widetilde
Q).$ For universal scheme  ${\rm Univ\,}^{P(m)}G(V, r)$ over
Hilbert scheme  ${\rm Hilb \,}^{P(m)}G(V,r)$ we have the fibred
diagram
\begin{equation*} \xymatrix{\widetilde \Sigma_i \ar[dd]
\ar[rd] \ar[rr]&& {\rm Univ\,}^{P(m)}G(V,r) \ar[dd]\\
& {\rm Univ\,}^{P(m)}G(V,r)|_{\mu(\widetilde Q)} \ar[dd] \ar[ur]&\\
B_i \ar[rd] \ar[rr]&& {\rm Hilb \,}^{P(m)}G(V,r)\\
& \mu(\widetilde Q) \ar[ur]}
\end{equation*}
Define schemes $B_{i Q}$ and $\widetilde \Sigma_{i Q}$ as fibred
products $B_{i Q}:=B_i \times_{\mu(\widetilde Q)} \widetilde Q$
and $\widetilde \Sigma_{i Q}:=\widetilde \Sigma_i
\times_{\widetilde Q} B_{i Q}.$ Let  $\widetilde \beta: \widetilde
\Sigma_{iQ} \to \widetilde \Sigma_i$ be the projection onto the
factor. Define $\widetilde {\mathbb E}_{iQ}:=\widetilde
\beta^{\ast}\widetilde {\mathbb E}_i$. Then there is a morphism
 $\Phi_i: \widetilde \Sigma_{iQ} \to
B_{iQ} \times S$ obtained by the fibred product from the morphism
$\Phi: \widetilde \Sigma_Q \to \widetilde Q \times S$. By Serre's
theorem and by choice of the invertible sheaf $\widetilde {\mathbb
L}$ one has an epimorphism $\pi^{\ast} \pi_{\ast} (\widetilde
{\mathbb E}_{iQ}\otimes \widetilde {\mathbb L}^m)\otimes
\widetilde {\mathbb L}^{-m} \twoheadrightarrow \widetilde {\mathbb
E}_{iQ}.$ Let $B_0$ be the open subset of the scheme $B$, where
 $\widetilde \Sigma$ is locally trivial and a fibre is isomorphic to
 $S$. Refining the cover  $\{B_i\}$ if necessary, we come to
 epimorphisms $V
\otimes L^{-m} \boxtimes {\mathcal O}_{B_{i0}}\twoheadrightarrow
\widetilde {\mathbb E}_{iQ}|_{B_{i0}}$ on overlaps $B_{i0}=B_0
\cap B_i$. Then there is a morphism $q_{i0}: B_{i0} \to {\rm
Quot\,}^{rp(m)}(V\otimes {\mathbb L}^{-m})$. Since $\widetilde
{\mathbb E}_{iQ}|_{B_{i0}}=\widetilde {\mathbb E}_{iQ}|_{B_{i0}}$
are families of semistable locally free sheaves, the morphism
$q_{i0}$ factors through the subscheme  $Q'$.

Form a closure  $\overline{q_{i0} B_{i0}}\supset q_{i0} B_{i0}$ in
the scheme  $Q'$, and a product $B_i \times \overline{q_{i0}
B_{i0}}$. Let  $\widetilde B_i$ be a closure of the subset of
points of the view $(b, q_{i0}(b)), b\in B_{i0}$ in the product
$B_i \times \overline{q_{i0} B_{i0}}$.  The projection $\rho:
\widetilde B_i \to B_i$ onto the first factor is a birational
morphism of integral schemes. Introduce the notation for the
composite map $q_i: \widetilde B_i \hookrightarrow B_i \times
\overline{q_{i0} B_{i0}}\stackrel{pr_2}{\longrightarrow}
\overline{q_{i0} B_{i0}} \subset Q'$.  Other necessary notations
are fixed by the following fibred diagram
\begin{equation*} \xymatrix{\widetilde \Sigma\!\!\!\Sigma_{iQ}
\ar[d]_{\widetilde \Phi_i} \ar[r]^{\widetilde \rho }& \widetilde \Sigma_{iQ} \ar[d]^{\Phi_i} \ar[r]& \widetilde \Sigma_Q \ar[d]^{\Phi}\\
\widetilde B_i \times S \ar[r]^{\rho \times id_S}& B_{iQ} \times
S\ar[r]& \widetilde Q \times S}
\end{equation*}
 Consider sheaves
${\mathbb E}_{iQ}:=(\widetilde \Phi_{i\ast} \widetilde
\rho^{\ast}\widetilde {\mathbb E}_{iQ})^{\vee \vee}$ and $(q_i,
id_S)^{\ast} {\mathbb E}_Q$. Both sheaves are reflexive on
integral and normal scheme $\widetilde B_i \times S$ and coincide
on open subset out of a subscheme of codimension 3 where they are
nonlocally free. Then ${\mathbb E}_{iQ}=(q_i, id_S)^{\ast}
{\mathbb E}_Q,$ what, in particular, proves that the sheaf
${\mathbb E}_{iQ}$ is flat over $\widetilde B_i$.

Composite maps $\widetilde B_i \to B_{iQ} \to \widetilde Q \to Q'$
factor through morphisms $B_i \to  \overline M$ under formation of
$PGL(V)$-quotient, since isomorphic semistable admissible pairs
correspond to $PGL(V)$-equivalent points. Morphisms  $B_i \to
\overline M$ are glued together in the morphism  $B \to \overline
M$. The dual morphism in the dual category $(Scheme\!s_k)^o$ leads
to the natural transformation $\omega:\underline{\overline M} \to
\underline F'.$

\end{document}